\documentclass[10pt,openbib]{article}
\usepackage{amsmath}
\usepackage{amsfonts}
\usepackage{amssymb}
\usepackage{graphicx}
\newcommand{\R}{\mathbb{R}} 
\newcommand{\C}{\mathbb{C}}
 
\newcommand{\CP}{\mathbb{C}\mathbb{P}}

\newcommand{\Z}{\mathbb{Z}}

\newcommand{\dee}{\mathop{\! \,\mathrm{d} \!}\nolimits}
\newcommand{\comp}{\raisebox{0pt}{$\scriptstyle\circ \, $}}
\newcommand{\setrule}{\, \rule[-4pt]{.5pt}{13pt}\, }

\newcommand{\onehalf}{\mbox{$\frac{\scriptstyle 1}{\scriptstyle 2}$}}

\newcommand{\ttfrac}[2]{\mbox{$\frac{{\scriptstyle #1}}{{\scriptstyle #2}}$}}
\newcommand{\circdot}{\raisebox{2pt}{\tiny$\, \bigodot$}}
\newcommand{\vvee}{\mbox{\tiny $\vee $}}

\allowdisplaybreaks

\begin{document}
\thispagestyle{empty}
\begin{center}
{\Large \textbf{The ${\mathbf{A}}_5$ Hamiltonian}}
\mbox{}\vspace{.05in} \\ 
\mbox{}\\
Richard Cushman\footnotemark 
\end{center}
\footnotetext{e-mail: r.h.cushman@gmail.com  \hfill printed: \today }

We discuss in detail the example of the holomorphic ${\mathbf{A}}_5$ Hamiltonian given in \cite{bates-cushman}. In particular, we use the finite group of discrete symmetries of the Hamiltonian to build a discrete subgroup of the $2$-dimensional Euclidean group whose orbit space is affine model of an affine Riemann surface, which is a 
level set of the ${\mathbf{A}}_5$ Hamiltonian. Our treatment differs from the one given in \cite{cushman} because the map, which develops the affine Riemann surface, is not constructed from a Schwarz-Christoffel transformation. Also we do not describe the 
monodromy of the ${\mathbf{A}}_5$ Hamiltonian, as it is fully discussed in 
\cite{bates-cushman}. %

\section{The Hamiltonian system}

Consider the holomorphic Hamiltonian system $(H, {\C }^2, \Omega = \dee \xi \wedge \dee \eta )$ 
with holomorphic Hamiltonian 
\begin{equation}
H: {\C }^2 \rightarrow \C : (\xi ,\eta ) \mapsto \onehalf {\eta }^2 + {\xi }^6 -1. 
\label{eq-one}
\end{equation}
The Hamiltonian vector field $X_H$ on $({\C }^2, \Omega )$ corresponding to $H$ is 
$\eta \frac{\partial }{\partial \xi } - 6{\xi }^5 \frac{\partial }{\partial \eta}$. A  
holomorphic integral curve 
$\gamma : \C \rightarrow {\C}^2: \tau \mapsto \big( \xi (\tau ), \eta (\tau ) \big) $ 
of $X_H$ satisfies the holomorphic form of Hamilton's equations 
\begin{equation}
\begin{array}{l}
{\displaystyle \frac{\dee \xi}{\dee \tau }  = \frac{\partial H}{\partial \eta } = \eta}  \\
\rule{0pt}{18pt} {\displaystyle \frac{\dee \eta }{\dee \tau }  = - \frac{\partial H}{\partial \xi } = -6 {\xi }^5}. 
\end{array}
\label{eq-onestar}
\end{equation}
 Because $H$ is constant on the holomorphic integral curves of $X_H$, $X_H$ is 
a holomorphic vector field on the affine Riemann surface
\begin{equation}
S = H^{-1}(0) = \{ (\xi ,\eta ) \in {\C }^{2} \setrule 
\onehalf {\eta }^2 = 1 - {\xi }^6 \} .
\label{eq-two}
\end{equation}

\section{Topology of $\boldsymbol{S}$}

To begin the study of the geometry of the integral curves of $X_H$, we \linebreak 
determine the topology of the affine Riemann surface $S \subseteq {\C }^2$ 
(\ref{eq-two}). The projective Riemann surface $\mathcal{S} \subseteq {\CP}^2$ corresponding to $S$ is defined by 
\begin{displaymath}
G(\xi ,\eta , \zeta ) = \onehalf {\eta }^2 {\zeta }^4 +{\xi }^6 - {\zeta }^6 = 0.
\end{displaymath}
The singular points of $\mathcal{S}$ are determined by 
\begin{displaymath}
\frac{\partial G}{\partial \xi } = 6{\xi }^5 = 0 ; \, \, 
\frac{\partial G}{\partial \eta } = {\zeta }^4\eta = 0 ; \, \, \mathrm{and} \, \, 
\frac{\partial G}{\partial \zeta } = 2{\zeta }^3({\eta }^2 - 3{\zeta }^2) = 0.
\end{displaymath}
Thus $\xi =0$. If $\zeta \ne 0$, then the second equation gives $\eta = 0$ and 
the third $0 = {\eta }^2 - 3{\zeta}^2 = -3{\zeta }^2$, that is, $\zeta =0$, which is 
a contradiction. So $\zeta =0$ and $\eta $ is arbitrary. Hence $[0:1:0]$ is the 
only singular point of $\mathcal{S}$. Thus $S = \mathcal{S} \cap \{ \zeta = 1 \} = 
\mathcal{S} \setminus \{ [0:1:0] \} $ is nonsingular.  

\par \noindent \hspace{1in}\begin{tabular}{l}
\includegraphics[width=220pt]{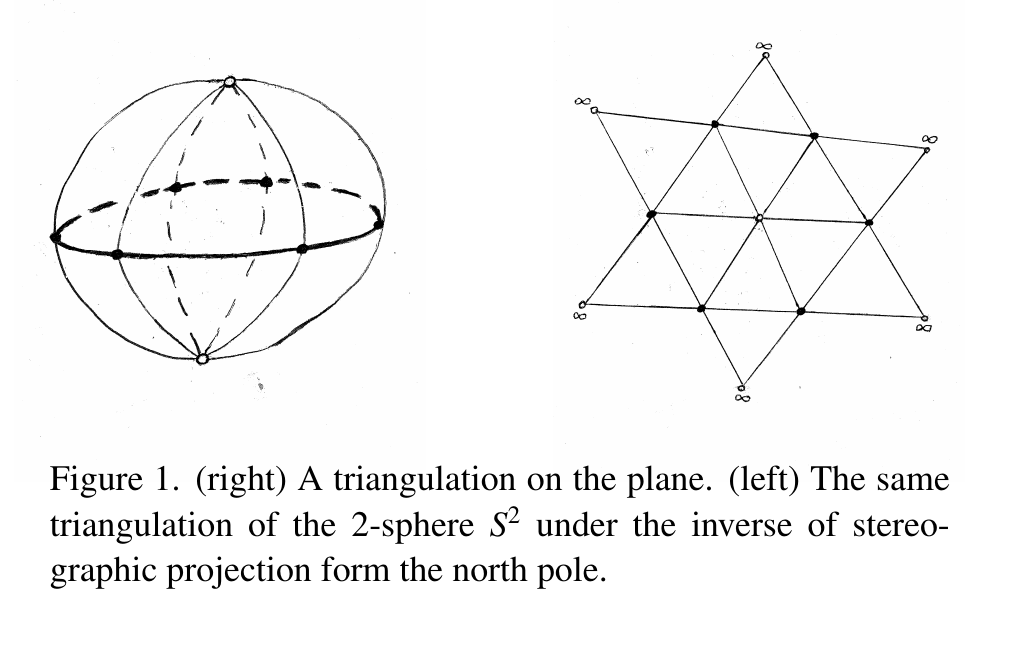}
\vspace{-.15in}
\end{tabular} \medskip 

\noindent \textbf{Lemma 2.1} The projection mapping 
\begin{displaymath}
\overline{\pi } : \CP^2 \rightarrow \CP: [\xi : \eta : \zeta ] \rightarrow [\xi : \zeta ],  
\end{displaymath}
when restricted to $\mathcal{S}$ presents $\mathcal{S}$ as a branched double covering of $\CP$ with branch points at $\{ [\xi : 0 : 1 ] \in \CP^2 \setrule \, 
{\xi }^6 = 1 \} $.  \medskip 

\noindent \textbf{Proof.} Observe that 
the affine map corresponding to $\overline{\pi}$ is the projection map  
\begin{equation}
\pi : \CP^2 \cap \{ \zeta = 1 \} = {\C }^2 \rightarrow \CP \cap \{ \zeta =1 \} = \C : 
 (\xi , \eta ) \mapsto \xi , 
 \label{eq-2zero}
\end{equation}
whose tangent at $(\xi , \eta )$ is $T_{(\xi ,\eta )}\pi (X \frac{\partial }{\partial \xi } + 
Y \frac{\partial }{\partial \eta }) = X \frac{\partial }{\partial \xi }$. 
The complex tangent space $T_{(\xi , \eta )}S$ to $S$ at $(\xi , \eta )$ is equal to 
$\ker \dee H(\xi , \eta )$, which is spanned by $X_H(\xi , \eta ) = 
\eta \frac{\partial }{\partial \xi } - 6{\xi }^5 \frac{\partial }{\partial \eta}$. Thus 
\begin{equation}
\big( T_{(\xi , \eta )}\pi |_S \big) (X_H(\xi , \eta )) = \eta \frac{\partial }{\partial \xi }. 
\label{eq-threestar}
\end{equation}
Hence $T_{(\xi ,\eta )}\pi |_S$ is a complex linear isomorphism of $T_{(\xi , \eta )}S$ 
onto $T_{\pi (\xi , \eta )}\C $ if and only if $(\xi , \eta ) \in S^{\dagger} = 
S \setminus \{ \eta = 0 \} = S \setminus \{ (\xi ,0) \in S \setrule \, {\xi }^6 =1 \}$. 
So $\pi |_{S^{\dagger}}$ is a local holomorphic diffeomorphism of $S^{\dagger}$ 
onto ${\C }^{\dagger } = \C \setminus \{ {\xi }^6 =1 \}$. It is a covering map having two sheets, since $({\pi }|_{S^{\dagger}})^{-1}(\xi ) = \{ (\xi , \pm \eta ) \in S^{\dagger} \} $ for every $\xi \in {\C}^{\dagger }$. The six points 
${\xi }_k = {\mathrm{e}}^{2\pi i k/6}$ for $k=0, \ldots , 5$ are 
branch values of the projection mapping ${\pi }_{|S}$ because near 
${\xi }_j$ we have
\begin{align*}
\eta & = (1-{\xi}^6)^{1/2} = \prod^5_{k=0}(\xi - {\xi }_k)^{1/2} \\
& = ({\xi }_j -{\xi }_k)^{1/2}( \xi -{\xi }_j)^{1/2}
\prod_{k\ne j} \big( 1 + ({\xi }_j- {\xi }_k)^{-1}(\xi - {\xi }_j) \big) ^{1/2}. 
\end{align*}
Each factor in the last product is a holomorphic function of $\xi -{\xi}_j$ when 
$|\xi -{\xi }_j|$ is sufficently small. \hfill $\square $ \medskip  %

The following argument shows that the genus of $\mathcal{S}$ is $2$. 
In figure 1 (right) the points at $\infty$ are identified and the edges from $\infty$ to a vertex of the regular hexagon are identified. With these identifications figure 1 (right) gives the triangulation of $S^2 = \CP$ in figure 1 (left) having 
V$\, = 8 $ vertices; E$\, = 18$ edges; F$\, =12$ faces. Thus the Euler characteristic $\chi $ of $S^2$ is $V-E+F = 8 -18 +12 =2$, which is $2 -2g$. Hence the genus $g$ of $S^2$ is $0$, as expected.  Taking two copies of figure 1 (right) with the same identifications as above and then identifying the darkened vertices of the regular hexagon gives a triangulation of $\mathcal{S}$ having $V^{\prime }\, = 10 = 2V -6 $ vertices; $E^{\prime}\, = 36 = 2E$ edges; $F^{\prime }\, = 24 =2F$ faces. Thus the Euler characteristic ${\chi }'$ of $\mathcal{S}$ is 
${\chi }^{\prime } = V^{\prime }-E^{\prime }+F^{\prime }  =-2$. Hence $-2 = 2-2g$, where $g$ is the genus of $\mathcal{S}$. So $g=2$. \medskip
 
Let $K$ be the closed stellated regular hexagon formed by repeatedly rotating the closed quadrilateral $Q' =OD'C\overline{D'}$ by $R$ through an angle $2\pi /6$, see figure 2. We now use the stellated regular hexagon $K$ to construct a Riemann surface of the same genus as $S^{\dagger}$, see 
\cite{richens-berry}.  \medskip 

Let $R$ be a rotation of $\C $ given by multiplication by ${\mathrm{e}}^{2\pi i/6}$ and 
let $U$ be the reflection given by complex conjugation. Let $G^{\vvee}$ be the group generated by the reflections 
$S_k = R^k SR^{-k} = R^{2k+1}U$ for $k=0,1, \ldots , 5$. Here $S = RU= S_0$ is the reflection, which leaves the closed ray $\ell  = \{ t {\mathrm{e}}^{2  \pi i /6} \setrule \, t \in OD'  \}$ fixed, and $S_k$ is the 
reflection, which leaves the ray $R^k\ell $ fixed. Note that $G^{\vvee} = \langle R, U \setrule \, 
R^6 = e = U^2, \, RU = UR^{-1} \rangle $. Define an equivalence relation on 
$K$ by saying that two points $x$ and $y$ in $K$ are \emph{equivalent}, $x \sim y$, 
if and only if 1) $x$ and $y$ lie on $\partial \, K$ with $x$ on 
the closed edge $E$ and $y$ on the closed edge of $K$, 
where $y = S_m(x) \in S_m(E)$ for some reflection 
$S_m \in G^{\vvee}$ or 2) if $x$ and $y$ lie in the interior of 
$K$ and $x =y$. Two edges of $K$ are equivalent if they contain equivalent points. Geometrically 
two edges are equivalent if they extend to lines in $\C$ which are parallel. 
Let $K^{\sim }$ be the space of equivalence classes of points in $K$ and let 
\begin{equation}
\rho : K \rightarrow K^{\sim } : p \mapsto [p]
\label{eq-s3threestar}
\end{equation}
be the identification map which sends a point $p \in K$ to its equivalence class $[p]$, which contains $p$. Give $K$ the topology induced from $\C $. Placing the quotient topology on 
$K^{\sim }$ turns it into a compact connected topological 
manifold without boundary, see the argument below. Let $K^{\ast }$ be $K$ less its vertices. \medskip

\noindent \textbf{Proposition 2.2} With $O$ the center of $K^{\ast }$ the identification space $(K^{\ast } \setminus  \mathrm{O} )^{\sim} =  \rho (K^{\ast } \setminus  \mathrm{O} )$ is a connected $2$-dimensional smooth manifold without boundary. \medskip 
 
\noindent \textbf{Proof.} Let $E_{+}$ be an open edge of $K^{\ast }$. For $p_{+} \in E_{+}$ let $D_{p_{+}}$ be a disk in $\C $ with center at $p_{+}$, which does not contain a vertex of $K$. Set $D^{+}_{p_{+}} = K^{\ast } \cap D_{p_{+}}$. Let $E_{-}$ be an open edge of $K^{\ast }$, which is equivalent to $E_{+}$ via the reflection $S_m$. Let $p_{-} = S_m(p_{+})$ and set $D^{-}_{p_{-}} = S_m(D^{+}_{p_{+}})$. Then $V_{[p]} =  
\rho (D^{+}_{p_{+}} \cup D^{-}_{p_{-}}) $ is an open neighborhood of $[p] = 
[p_{+}] = [p_{-}]$ in $(K^{\ast } \setminus  \mathrm{O} )^{\sim }$, which is a smooth $2$-disk, since the identification mapping 
$ \rho $ is the identity on $\mathrm{int}\, K^{\ast }$. It follows that 
$(K^{\ast } \setminus  \mathrm{O}  )^{\sim }$ is a smooth $2$-dimensional manifold without boundary. \hfill $\square $ \medskip %

\noindent \textbf{Proposition 2.3} The identification space $K^{\sim }$ is a connected compact 
topological manifold with a conical singular point at each vertex of $K$ and at the center $\mathrm{O}$ of $K$. \medskip 

\noindent \textbf{Proof.} We now handle the vertices of $K$. Let $v_{+}$ be a vertex of $K^{\ast }$ and set 
$D_{v_{+}} =\widetilde{D} \cap K^{\ast }$, where $\widetilde{D}$ is a disk in $\C $ 
with center at the vertex $v_{+}= r_0{\mathrm{e}}^{i \pi {\theta }_0}$. The map 
\begin{displaymath}
W_{v_{+}}: D_{+} \subseteq \C \rightarrow D_{v_{+}} \subseteq \C: 
r{\mathrm{e}}^{i \pi \theta } \mapsto |r - r_0| {\mathrm{e}}^{i \pi s (\theta - {\theta }_0)} 
\end{displaymath}
with $r \ge 0$ and $0 \le \theta \le 1$ is a homeomorphism, which sends the wedge with angle $\pi $ to the wedge with angle $\pi s$. The latter wedge is formed by the closed edges $E'_{+}$ and $E_{+}$ of $K^{\ast }$, which are adjacent at the vertex $v_{+}$ such that ${\mathrm{e}}^{i\pi s}E'_{+} = E_{+}$ with the edge $E'_{+}$ being swept out through 
$\mathrm{int}\, K^{\ast }$ during its rotation to the edge $E_{+}$. Because $K^{\ast }$ is a rational stellated regular hexagon, the value of $s$ is a rational number for each vertex of 
$K^{\ast }$. Let $E_{-} = S_m(E_{+})$ be 
an edge of $K^{\ast }$, which is equivalent to $E_{+}$ and set 
$v_{-} = S_m(v_{+})$. Then $v_{-}$ is a vertex of $K^{\ast }$, which is 
the center of the disk $D_{v_{-}} = S_m(D_{v_{+}})$. Set $D_{-} = 
{\overline{D}}_{+}$. Then $D = D_{+} \cup D_{-}$ is a disk in $\C $. The map 
$W: D \rightarrow  \rho (D_{v_{+}} \cup D_{v_{-}} )$, where $W|_{D_{+}} =  
\rho \comp W_{v_{+}}$ and $W|_{D_{-}} =  \rho \comp S_m \comp W_{v_{+}} \comp 
{\mbox{}}^{\raisebox{2pt}{$\overline{\rule{5pt}{0pt}}$} }$, is a homeomorphism of $D$ into a neighborhood 
$ \rho ( D_{v_{+}} \cup D_{v_{-}})$ of $[v] = [v_{+}] = [v_{-}]$ in 
$(K^{\ast })^{\sim}$. \medskip 

The center $\mathrm{O}$ of $K$ is a conical singular point of $K^{\sim}$ because it is a fixed point of the 
linear action of $G^{\vvee}$ on $\C$. \medskip 

Consequently, the identification space 
$(K^{\ast })^{\sim }$ is a topological manifold with a conical singularity at a point corresponding to a vertex or 
the center $\mathrm{O}$ of $K$. \hfill $\square $ \medskip

Let $\widetilde{G}= \langle R \setrule \, R^6 = e \rangle $. 
Then $\widetilde{G}$ is the abelian group $\Z  \bmod 6$. 
The usual $\widetilde{G}$-action 
$\widetilde{G} \times K \subseteq \widetilde{G} \times \C \rightarrow K \subseteq \C :(g,z) \mapsto g(z)$ 
preserves equivalent edges of $K$ and is free on 
$K^{\ast } \setminus  \mathrm{O} $. Hence it induces a $\widetilde{G}$ action on 
$(K^{\ast } \setminus  \mathrm{O}  )^{\sim }$, which is free and proper with 
orbit map  
\begin{displaymath}
\sigma :  (K^{\ast } \setminus  \mathrm{O}  )^{\sim } \rightarrow 
(K^{\ast } \setminus  \mathrm{O}  )^{\sim }/\widetilde{G} = 
{\widetilde{S}}^{\dagger}: z \mapsto z\widetilde{G}
\end{displaymath} 

\noindent \textbf{Proposition 2.4} The orbit space 
${\widetilde{S}}^{\dagger} = \sigma ((K^{\ast } \setminus  \mathrm{O}  )^{\sim})$ is a connected smooth complex manifold. \medskip  

\noindent \textbf{Proof.} Since the action of $\widetilde{G}$ on 
$(K^{\ast} \setminus  \mathrm{O} )^{\sim}$ is free and proper, its orbit 
space $\sigma ((K^{\ast } \setminus  \mathrm{O}  )^{\sim})$ is a smooth 
$2$-dimensional manifold that is the orientated. This orientation is induced from an orientation of $K^{\ast } \setminus  \mathrm{O}  $, which comes from the 
orientation $\C $. So ${\widetilde{S}}^{\dagger}$ has a complex structure, since each element of $\widetilde{G}$ is a conformal mapping of $\C $ into itself. 
\hfill $\square $ \medskip

\vspace{-.15in} \par \noindent \hspace{.2in}\begin{tabular}{l}
\includegraphics[width=300pt]{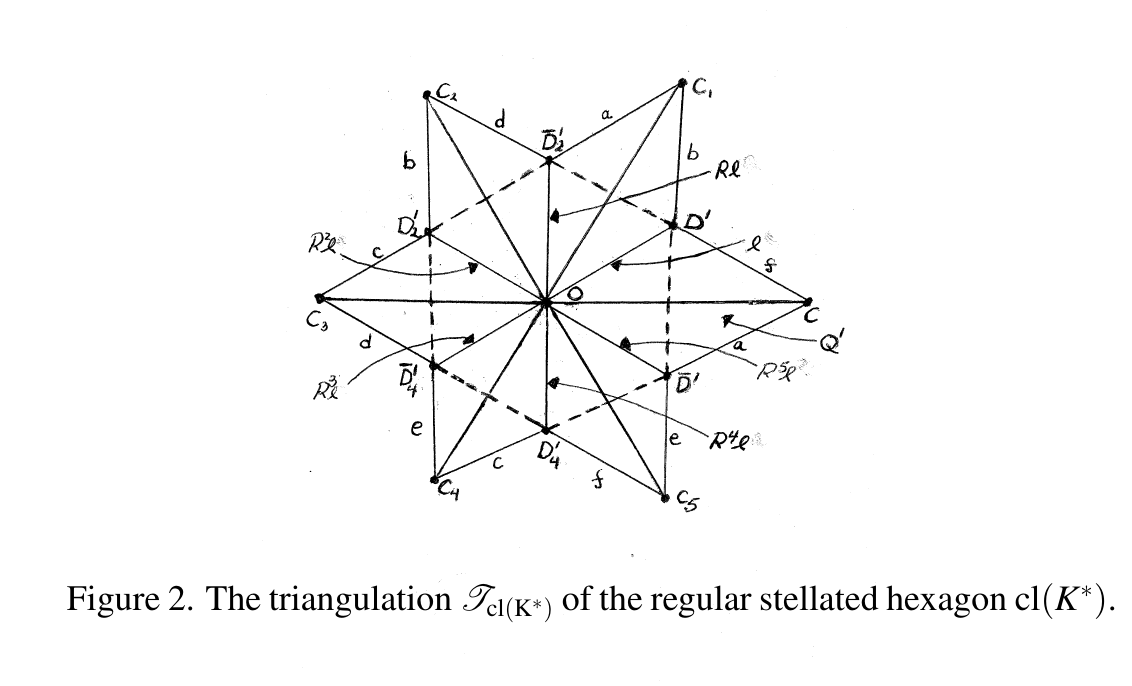}
\vspace{-.1in}
\end{tabular}  

Next we specify the topology of ${\widetilde{S}}^{\dagger}$. The stellated regular hexagon $K^{\ast } \setminus  \mathrm{O} $ less the origin has 
a triangulation ${\mathcal{T}}_{K^{\ast } \setminus  \mathrm{O} }$ made up of 
$12$ open triangles $R^j(\bigtriangleup OCD')$ and $R^j(\bigtriangleup OC{\overline{D}}')$ for $j =0,1, \ldots ,5$; $24$ open edges $R^j(OC)$, $R^j(O{\overline{D}}')$, $R^j(C{\overline{D}}')$, and $R^j(CD')$ for $j =0,1, \ldots , 5$; and $12$ vertices $R^j(D')$ and $R^j(C)$ for $j = 0,1, \ldots , 5$, see figure 2. \medskip 

Consider the set $\mathcal{E}$ of unordered pairs of equivalent closed edges of $K$, that is, $\mathcal{E}$ is the 
set $[E, S_k(E)]$ for $k=0,1, \ldots , 5$, where $E$ is a closed edge of $K$ and $S_k \in G^{\vvee}$. Table 1 lists the elements of $\mathcal{E}$. $G^{\vvee}$ acts on $\mathcal{E}$, namely, 
$g \cdot [E, S_k(E)] = [g(E), gS_kg^{-1}\big( g(E) \big) ] $, for $g \in G^{\vvee}$. 
This $G^{\vvee}$ action is well defined, since $gS_kg^{-1}$ is a reflection in the line $g(R^k\ell )$. This argument shows that $\widetilde{G}$ acts on $\mathcal{E}$, since $\widetilde{G}$ is a subgroup of $G^{\vvee}$.  \medskip 

\par \noindent \hspace{.3in} \begin{tabular}{lcl}
$a = \big[ \overline{D'}C, S_0( \overline{D'}C) = \overline{D'_2}C_1 \big] $ 
&\quad & 
$b = \big[ D'C_1, S_1(D'C_1) = D'_2C_2 \big] $ \\
\rule{0pt}{12pt}$d = \big[ \overline{D'_2}C_2, S_2(\overline{D'_2}C_2) = 
\overline{D'_4}C_3 \big] $ 
& \quad & 
$c = \big[ D'_2C_3, S_3(D'_2C_3) = D'_4C_4 \big] $ \\ 
\rule{0pt}{12pt}$e = \big[ \overline{D'_4}C_4, S_4(\overline{D'_4}C_4) = 
\overline{D'}C_5 \big] $ 
& \quad & 
$f = \big[ D'_4C_5, S_5(D'_4C_5) = D'C \big] $
\end{tabular} \medskip 

\par \noindent \hspace{.4in} \footnotesize{Table 1.} \parbox[t]{3.3in}{\footnotesize Elements of the set $\mathcal{E}$ of unordered pairs of equivalent closed edges of $K$. Here $D'_k =R^k(D')$ and $\overline{D'_k}= R^k(\overline{D'} )$ for $k=0,2,4$ and $C_k = R^k(C)$ for $k= \{ 0,1, \ldots , 5 \} $, see figure 2.}  \medskip 

\normalsize We now look at the $G^{\vvee}$-orbits on $\mathcal{E}$. We compute the $G^{\vvee}$-orbit of $d \in \mathcal{E}$ as follows.   
\begin{align}
(UR) \cdot d & = 
\big[ UR(\overline{D'_2}C_2), UR(S_2(\overline{D'_2}C_2) ) \big] = 
\big[ UR(\overline{D'_2}C_2), UR(\overline{D'_4}C_3)) \big] \notag 
\\
& = \big[ U(D'_2 C_3), U(D'_4C_4) \big] = 
\big[ \overline{D'_4}C_3, \overline{D'_2}C_2 \big] = d, \notag \\ 
R^2 \cdot d & = R^2 \cdot \big[ \overline{D'_2}C_2, S_2(\overline{D'_2}C_2) \big]  = \big[ R^2(\overline{D'_2}C_2), R^2S_2R^{-2}(R^2(\overline{D'_2}C_2) ) \big] \notag 
\\
& = \big[ \overline{D'_4}C_4, S_4( \overline{D'_4} C_4) \big] = 
\big[ \overline{D'_4}C_4, \overline{D'}C_5 \big] = e, \notag
\end{align}
and 
\begin{align}
R^4 \cdot d & =  
\big[ R^4(\overline{D'_4}C_2), R^4S_2R^{-4}(R^4(\overline{D'_2}C_2) )\big] \notag \\
&= \big[ \overline{D'}C , S_6(\overline{D'}C) \big] = 
\big[ \overline{D'}C, S_0(\overline{D'}C) \big] 
= \big[ \overline{D'}C, \overline{D'_2}C_1 \big] = a . \notag
\end{align}
Thus the $G^{\vvee}$ orbit $G^{\vvee} \cdot d$ of $d \in 
\mathcal{E}$ is $H \cdot d = \{ a, d, e \} $, which is a $\widetilde{G}$ orbit, since $H = G^{\vvee}/\langle UR \, | \, (UR)^2  = e \rangle = \langle V = R^2 \setrule \, V^3 = e  \rangle $ is a subgroup of $\widetilde{G}$. 
Similarly, the $G^{\vvee}$-orbit $G^{\vvee} \cdot f$ of $f \in \mathcal{E}$ is $H \cdot f = \{ b, c, f \} $. Since 
$G^{\vvee} \cdot d \, \cup \, G^{\vvee} \cdot f = \mathcal{E}$, we have found all 
$G^{\vvee}$-orbits and hence all the $\widetilde{G}$-orbits on $\mathcal{E}$.  \medskip 

The end points of the elements of the orbit $G^{\vvee}\cdot d$ are $a = 
\{ {\overline{D}}' , C, {\overline{D}}'_2, C_1 \} ; 
d = \{ {\overline{D}}'_2, C_2, {\overline{D}}'_4, C_3 \}$; and $e = \{ {\overline{D}}'_4, C_4 ; {\overline{D}}', C_5 \} $, 
see figure 2. Thus $\{ {\overline{D}}', {\overline{D}}'_2, {\overline{D}}'_4 \} ,$ $ \{ C, C_2, C_4 \}, 
\{ {\overline{D}}'_2, {\overline{D}}'_4, {\overline{D}}' \} $ and $\{ C_1, C_3, C_5 \}$ are $G^{\vvee}$ orbits of 
vertices of K. They are also $\widetilde{G}$ orbits. Similarly the end points of the orbit $G^{\vvee}\cdot f$ are 
$b = \{ D', C_1, D'_3,$ $C_2 \}; c = \{ D'_3, C_3 , D'_4, C_4\} $ and $f = \{ D'_4, C_5 , D' , C \} $. So 
$\{ D' , D'_3, D'_4 \}$, $\{ C_1,C_3,$ $ C_5 \}$, $\{ D'_3, D'_4, D' \}$ and $\{ C_2,C_4,C \}$ are $G^{\vvee}$ orbits 
of vertices of $K$. They are also $\widetilde{G}$ orbits. \medskip  

To determine the topology of the $\widetilde{G}$ orbit space 
${\widetilde{S}}^{\dagger}$ we find a triangulation of its closure. Note that the triangulation ${\mathcal{T}}_{K^{\ast } \setminus  \mathrm{O} }$ of 
$K^{\ast } \setminus  \mathrm{O} $, illustrated in figure 2, is $\widetilde{G}$-invariant. Its image under the identification map $ \rho $ is a $\widetilde{G}$-invariant triangulation 
${\mathcal{T}}_{(K^{\ast } \setminus  \mathrm{O})^{\sim} }$ of 
$(K^{\ast } \setminus  \mathrm{O} )^{\sim }$ with vertices 
$\rho (v)$, where $v$ is a vertex of $K$; open edges 
$\rho (E)$ having $\rho (O)$ as an end point, where $E$ is an edge of 
${\mathcal{T}}_{K^{\ast } \setminus O}$ having $O$ as an end point; 
or $\rho ([F,F'])$, where $[F,F']$ is an unordered pair of equivalent edges in $K$; 
open triangles $\rho (T)$, where $T$ is a triangle in 
${\mathcal{T}}_{K^{\ast } \setminus O}$. The triangulation 
${\mathcal{T}}_{(K^{\ast }\setminus O)}$ is invariant under the induced 
$\widetilde{G}$ action on $(K^{\ast } \setminus O)^{\sim}$. It follows that 
$\sigma ( \rho (v))$, $\sigma ( \rho (E))$ or $\sigma ( \rho ([F,F']))$, and 
$\sigma ( \rho (T))$ is a vertex, an open edge, and an open triangle, respectively, of a triangulation 
${\mathcal{T}}_{{\widetilde{S}}^{\dagger}} = 
\sigma ({\mathcal{T}}_{(K^{\ast } \setminus  \mathrm{O} )^{\sim }})$ of 
$\widetilde{S}$. The triangulation 
${\mathcal{T}}_{{\widetilde{S}}^{\dagger}}$ has 1) $4$ vertices, 
corresponding to the $\sigma \comp \rho $ image of the $\widetilde{G}$ orbits 
$\{ D' , D'_3, D'_4 \}$, $\{ C_1,C_3, C_5 \}$, $\{ \overline{D}'_2, \overline{D}'_4, \overline{D}' \}$ and 
$\{ C_2,C_4,C \}$ of vertices of $K$. 2) $8$ open edges corresponding to $\sigma \comp \rho $ image of the 
$\widetilde{G}$ orbits of the two 
edge pairs $d$ and $f$ of $K$ and $\sigma \comp \rho $ image of the six $\widetilde{G}$ orbits $\{ R^{2j}(OC)\} $ 
and $\{R^{2j}(OD')\} $ for $j=0,1,2$ and 
$\{ R^{2j}(R(OC)) \} $ and $\{ R^{2j}(R(OD')) \} $ for $j=0,1,2$ and $\{ R^{2j}(CD') \}$ and 
$\{ R^{2j}(C_1\overline{D}'_2)\} $ for $j=0,1,2$;\footnote{Note that the end points of the edges of these six 
$\widetilde{G}$ orbits are $\widetilde{G}$ orbits of end points of equivalent edges or the center $\mathrm{O}$ of 
$K$.} 3) $2$ open triangles corresponding to the $\sigma \comp \rho $ image of $\widetilde{G}$ orbits of the triangles $\bigtriangleup OCD'$ and $\bigtriangleup OC\overline{D'}$ of $K$. Thus the Euler characteristic $\chi (\widetilde{\mathcal{S}}')$ of
$\widetilde{\mathcal{S}}'$ is $4 - 8 + 2 = -2$. Since ${\widetilde{S}}^{\dagger}$ is a $2$-dimensional smooth real manifold, $\chi ({\widetilde{S}}^{\dagger}) = 2 - 2g$, where $g$ is the genus of ${\widetilde{S}}^{\dagger}$. Hence $g=2$. So 
${\widetilde{S}}^{\dagger}$ is a smooth $2$-sphere with $2$ handles,  less a finite number of points, which lies in a compact topological space $K^{\sim }/\widetilde{G}$, that is its closure.  

\section{The developing mapping}

Consider the mapping 
\begin{equation}
F : {\C }^{\dagger } = \C \setminus \{ {\xi }^6 = 1 \} \rightarrow \C : \xi \mapsto z = 
\frac{1}{\sqrt{2}} \int^{\xi }_0 \frac{\dee w }{\sqrt{1- w^6}}. 
\label{eq-four}
\end{equation}
$F$ is holomorphic except at the sixth roots of unity $\{ \xi \in \C \setrule \, 
{\xi }^6 =1 \}$, since $\dee z = \frac{1}{\sqrt{1- {\xi }^6}} \dee \xi $. In addition, 
$F$ is a local holomorphic diffeomorphism on ${\C }^{\dagger}$, because $F' \ne 0$. Let $\omega = {\mathrm{e}}^{2\pi \mathrm{i}/6}$ and set $R: \C \rightarrow \C : \xi \mapsto \omega \, \xi $. Since 
\begin{align*}
F(\omega \xi ) & = \int^{\omega \xi}_0 \frac{\dee w}{\sqrt{1-w^6}} = 
\int^{\xi}_0 \frac{\dee w'}{\sqrt{1-(w')^6}}, \, \, \mbox{where $w = \omega w'$} \\
& = {\omega }^{-1}\, \int^{\xi}_0 \frac{\dee w}{\sqrt{1-w^6}} = {\omega }^{-1} F(\xi ), 
\end{align*}
the map $F$ intertwines the ${\Z }_6$ action on ${\C}^{\dagger} =\C \setminus \{ {\xi }^6 =1 \}$ generated by $R$, that is, $F(R\xi ) = R^{-1}(F(\xi))$ for every $\xi \in {\C }^{\dagger}$. So $F({\omega }^k) = {\omega }^{-k}F(1)$ for $k=0,1, \ldots , 5$. Let $T$ be the closed triangle with edges $\overline{0\, 1}$, 
$\overline{1\, \omega }$, and $\overline{0\, \omega }$, see figure 3 (left).  

\par \noindent \hspace{.75in}\begin{tabular}{l}
\includegraphics[width=250pt]{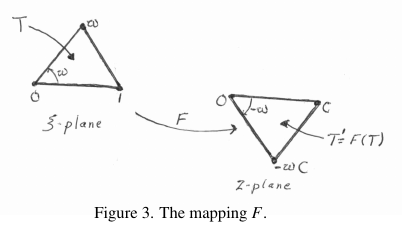}
\vspace{-.1in}
\end{tabular}  
\par \noindent The image of $T$ under the mapping $F$ is 
the triangle $T'$ with edges $F(\overline{0\, 1}) = \overline{0\, C}$, $F(\overline{\omega \, 1}) = 
\overline{C\, (-\omega )C}$, and $F(\overline{0\, \omega}) = {\omega }^{-1}F(\overline{0 \, 1}) = \overline{0(-\omega C)}$, see figure 3 (right). Here $C = F(1) = \int^1_0\frac{\dee w}{\sqrt{1-w^6}}$. \medskip

Let $H = {\bigcup}^5_{\nu =0} R^{\nu }(T)$. Then $H$ is a closed regular hexagon whose vertices are at the sixth roots of unity. Since $H$ is simply connected, the complex square root is single valued there. Hence 
$\mathcal{D} = ({\pi }_{|S^{\dagger}})^{-1}(H)$ is one sheet of the 
projection map ${\pi }_{|S^{\dagger}}$ (\ref{eq-2zero}). Since ${\pi }_{|S^{\dagger}}$ is a twofold branched covering map, the other sheet is $\mathcal{D}' = ({\pi }_{|S^{\dagger}})^{-1}(K \setminus H)$, where 
$K$ is the stellated regular hexagon formed by adding an equilateral 
triangle with an edge in common with the hexagon $H$ for each edge of 
$H$, see figure 1 (right). The image of ${\C}^{\dagger}$ under the mapping $F$ 
(\ref{eq-four}) is the closed stellated regular hexagon $K$ in $\C$ with center at the origin and side length $C$. \medskip 

Using (\ref{eq-threestar}) we get 
\begin{align}
(T_{\xi }F \comp \pi |_{S^{\dagger } })(X_H(\xi , \eta )) &  = 
T_{\xi }F \big( T_{(\xi , \eta)}\pi |_{S^{\dagger } }(X_H(\xi , \eta )) \big) 
= T_{\xi }F (\eta \frac{\partial }{\partial \xi }) = 
\frac{\partial }{\partial z}, 
\label{eq-fivestar}
\end{align}
since equation (\ref{eq-four}) gives $\dee z = \dee F = \frac{1}{\eta }\, \dee \xi $, where $\eta = \sqrt{2(1-{\xi }^6)}$. Thus the mapping 
\begin{equation}
\delta = F \comp {\pi }_{|S^{\dagger}}: S^{\dagger} \subseteq {\C}^2 
\rightarrow K \subseteq 
{\C }^{\vvee} = \C \setminus \{ z^6 = C^6\} \subseteq \C 
\label{eq-sixstar}
\end{equation}
\emph{straightens} the holomorphic vector field $X_H$ on $S^{\dagger}$, that is, on $\mathcal{D}$ 
and ${\mathcal{D}}'$. Moreover, ${\delta }_{|\mathcal{D}}$ and 
${\delta }_{|{\mathcal{D}}'}$ are holomorphic diffeomorphisms of $\mathcal{D}$ and 
${\mathcal{D}}'$ onto the hexagon $H$ less its vertices and $K \setminus H$, less its vertices, respectively. \medskip

Consider the hermitian metric $\gamma = \dee z \circdot \dee \overline{z}$ on 
${\C}^{\vvee}$. Pulling ${\gamma }_{|K}$ back by the map 
$\delta $ (\ref{eq-sixstar}) gives the hermitian metric $\Gamma = 
\frac{1}{\eta } \dee \xi \circdot \frac{1}{\overline{\eta }} \dee \overline{\xi }$ on 
$S^{\dagger}$, since 
\begin{align}
\Gamma (X_H(\xi , \eta ), X_H(\xi , \eta )) & =  
\frac{1}{\eta } \dee \xi \big( \eta \frac{\partial }{\partial \xi } -6{\xi }^5 
\frac{\partial }{\partial \eta }\big) \circdot \frac{1}{\overline{\eta }} 
\dee \overline{\xi }\big( \overline{\eta \frac{\partial }{\partial \xi } -6{\xi }^5 
\frac{\partial }{\partial \eta }} \big) \notag \\
&\hspace{-.75in} = \frac{1}{\eta } \dee \xi \big( \eta \frac{\partial }{\partial \xi } \big) \circdot 
\frac{1}{\overline{\eta }} \dee \overline{\xi } \big( \overline{\eta \frac{\partial }{\partial \xi }}\big) 
= \dee {\xi }( \frac{\partial }{\partial \xi }) \circdot \dee \overline{\xi }(\frac{\partial }{\partial \overline{\xi }}) = 1. \notag  
\end{align}
The metric $\Gamma $ on $S^{\dagger}$ is flat, because the metric ${\gamma }_{|K}$ on 
${\C }^{\vvee}$ is flat. The map $\delta $ (\ref{eq-sixstar}) is a holomorphic isometry and is also a local holomorphic diffeomorphism. Hence $\delta $ is a \emph{developing map}. Since the real integral curves of $\alpha \frac{\partial }{\partial z}$, $\alpha \in S^1$, $\alpha \in S^1$ on $({\C }^{\vvee}, \gamma )$ are geodesics, which are straight lines that make an angle $\theta $, where $\alpha = {\mathrm{e}}^{\mathrm{i}\theta }$, with the $x$-axis, the real integral curves of $X_{\alpha H}$ on $(S^{\dagger}, \Gamma )$ are geodesics for the hermitian metric $\Gamma $. Because the image under $(T\delta )^{-1}$ of the vector field $\alpha \frac{\partial }{\partial z}$ is the vector field $X_{\alpha H}$ on 
$S^{\dagger}$, the image of a real integral curve of $\alpha \frac{\partial }{\partial z}$ is the image of a real integral curve of $X_H$ and hence is a geodesic on $(S^{\dagger}, \Gamma )$. 

\noindent \section{Metric geometry of $\boldsymbol{S^{\dagger}}$}

In this section we discuss some properties of geodesics on the smooth affine 
Riemann surface $S^{\dagger}$. \medskip 

\noindent \textbf{Lemma 4.1} The mappings $\mathcal{R}: S^{\dagger} \subseteq 
{\C}^2 \rightarrow S^{\dagger} \subseteq {\C }^2: (\xi ,\eta ) \mapsto 
(\omega \xi , \eta )$, where $\omega = {\mathrm{e}}^{2\pi i/6}$ and 
$\mathcal{U}: S^{\dagger} \subseteq {\C }^2 \rightarrow S^{\dagger} \subseteq 
{\C }^2: (\xi , \eta ) \mapsto (\overline{\xi}, \overline{\eta })$ 
are isometries of $(S^{\dagger} , \Gamma )$. \medskip 

\noindent \textbf{Proof.} We compute. 
\begin{align}
{\mathcal{R}}^{\ast }\Gamma (\xi ,\eta )& = {\mathcal{R}}^{\ast}\big( \frac{1}{\eta }{\dee \xi} \circdot  
\frac{1}{\overline{\eta }}{\dee \overline{\xi}} \big) = 
\frac{1}{\eta }\dee \, (\omega \xi ) \circdot  
\frac{1}{\overline{\eta }}\dee \, \overline{(\omega  \xi )} 
= \frac{1}{\eta }{\dee \xi} \circdot  \frac{1}{\overline{\eta }}{\dee \overline{\xi}} = 
\Gamma (\xi , \eta ) . \notag 
\end{align}
$\mathcal{U}$ maps $H^{-1}(0)$ into itself, for if $(\xi , \eta ) \in H^{-1}(0)$, then 
$\onehalf {\overline{\eta }}^{\, 2} + {\overline{\xi }}^{\raisebox{-3pt}{\hspace{.5pt}${\scriptstyle 6 }$}} -1 = 
\overline{\onehalf {\eta }^2 +{\xi }^6 -1}$ $= 0$. 
So $(\overline{\xi }, \overline{\eta }) \in H^{-1}(0)$. The set $\{ {\xi }^6 = 1\} $ is mapped onto 
itself by $\mathcal{U}$. Hence $\mathcal{U}$ preserves $S^{\dagger}$. Now 
\begin{displaymath}
{\mathcal{U}}^{\ast }\Gamma (\xi , \eta ) = {\mathcal{U}}^{\ast }\big( \frac{1}{\eta }{\dee \xi} \circdot  
\frac{1}{\overline{\eta }}{\dee \overline{\xi}} \big) = 
\frac{1}{\overline{\eta }}{\dee \overline{\xi}} \circdot  
\frac{1}{\overline{\overline{\eta }}}{\dee \overline{\overline{\xi}}} = 
\frac{1}{\eta }{\dee \xi} \circdot  \frac{1}{\overline{\eta }}{\dee \overline{\xi}} = \Gamma (\xi ,\eta ) .
\end{displaymath}%
Thus $\mathcal{U}$ is an isometry of $(S^{\dagger}, \Gamma )$. \hfill $\square $ \medskip 

\noindent So the group $\widehat{\mathcal{G}} = \langle \mathcal{R}, \mathcal{U} \setrule \, {\mathcal{R}}^6 = {\mathcal{U}}^2 =e, \, \, \& \, \, \mathcal{R}\mathcal{U} = \mathcal{U}{\mathcal{R}}^{-1} \rangle $ is a group of isometries of $(S^{\dagger}, \Gamma )$. \medskip

\noindent \textbf{Lemma 4.2} The image of a geodesic on $(S^{\dagger} , \Gamma )$ under the action of the group of isometries $\widehat{\mathcal{G}}$ is a geodesic on $S^{\dagger}$. \medskip 

\noindent \textbf{Proof.} This follows because a geodesic is locally 
length minimizing, which is a property preserved by an isometry. 
\hfill $\square $ \medskip 

Let $\widehat{G}$ be the group of invertible linear maps of $\C $ into itself generated by 
\begin{displaymath}
R: \C \rightarrow \C : z \mapsto \omega z \, \, \, \mathrm{and} \, \, \, 
U: \C \rightarrow \C : z \mapsto \overline{z}, 
\end{displaymath}
which satisfy the relations $R^6 = U^2 = e$ and $RU = UR^{-1}$. The elements of 
$\widehat{G}$ preserve ${\C }^{\vvee} = \C \setminus \{ z^6 = C^6 \}$. \medskip 

\noindent \textbf{Claim 4.3} The developing map 
$\delta : S^{\dagger} \subseteq {\C }^2 \rightarrow {\C }^{\vvee}$ 
(\ref{eq-sixstar}) intertwines the $\widehat{\mathcal{G}}$ action $\Phi $ on $S^{\dagger}$ with the 
$\widehat{G}$ action $\phi $ on ${\C }^{\vvee}$. Specifically, for every $g \in \widehat{\mathcal{G}}$ and every $(\xi , \eta ) \in S^{\dagger}$ we have 
\begin{equation}
\delta \big( {\Phi }_g(\xi ,\eta ) \big) = {\phi }_{\psi (g)}\big( \delta (\xi ,\eta ) \big) ,  
\label{eq-seven}
\end{equation}
where $\psi : \widehat{\mathcal{G}} \rightarrow \widehat{G}$ is the isomorphism defined by $\psi (\mathcal{R}) = R$ and 
$\psi (\mathcal{U}) = U$. \medskip 

\noindent \textbf{Proof.} The following computation shows that equation (\ref{eq-seven}) holds when $g$ is $\mathcal{R}$ and $\mathcal{U}$, respectively. \medskip 

\noindent 1. $g = \mathcal{R}$. 
\begin{align}
\delta \big( {\Phi }_{\mathcal{R}}(\xi ,\eta ) \big) & = 
F(\omega \xi )  = \omega  F(\xi ) 
= R\big( \delta (\xi ,\eta ) \big) = {\phi }_{\psi (\mathcal{R})}\big( \delta (\xi ,\eta ) \big). \notag
\end{align}

\noindent 2. $g = \mathcal{U}$. 
\begin{align}
\delta \big( {\Phi }_{\mathcal{U}}(\xi , \eta ) \big) & = \delta (\overline{\xi }, \overline{\eta }) 
=F \big( \pi (\overline{\xi }, \overline{\eta }) \big) = F(\overline{\xi }) 
= \overline{F(\xi)} = 
U\big(\delta (\xi , \eta ) \big) = {\phi }_{\psi (\mathcal{U})}\big( \delta (\xi , \eta ) \big). 
\notag 
\end{align}
The fourth equality above follows by changing the variable $w$ to $\overline{w}$ 
in the integral $\frac{1}{\sqrt{2}}\int^{\xi }_0\frac{\dee w}{\sqrt{1-w^6}}$ \hfill 
$\square $. \medskip

\noindent \textbf{Corollary 4.3A} The developing map $\delta $ (\ref{eq-sixstar}) intertwines the action of 
$\widehat{\mathcal{G}}$ on $S^{\dagger}$ with the action $\phi $ of 
$\widehat{G}$ restricted to the stellated hexagon $K$. \medskip 

\noindent \textbf{Corollary 4.3B} $(S^{\dagger}, \Gamma )$ is geodesically incomplete. \medskip 

\noindent \textbf{Proof.} Since $(S^{\dagger }, \Gamma )$ is isometric to 
$({\C }^{\vvee}, \dee z \circdot \dee \overline{z} )$, which is 
$( ({\R }^2)^{\vvee},  \varepsilon = \dee x \circdot \dee x + \dee y \circdot \dee y )$, where 
\begin{displaymath}
{(\R }^2)^{\vvee } = {\R }^2 \setminus \{ (C\cos \ttfrac{2\pi k}{6}, C\sin \ttfrac{2\pi k}{6}) \in {\R }^2 \setrule 
\, k=0,1, \ldots , 5 \} , 
\end{displaymath}
it suffices to show that there are integral curves of the geodesic vector field 
$\frac{\partial }{\partial x}$ on $(({\R }^2)^{\vvee}, \varepsilon )$ 
which run off $({\R}^2)^{\vvee }$ in finite time. Consider the horizontal line segment 
\begin{displaymath}
\begin{array}{l}
\gamma : [0,1] \rightarrow {\R }^2: t \mapsto (1-t)(C\cos \frac{4\pi}{6}, 
C\sin \frac{4\pi}{6} ) +t (C \cos \frac{5\pi }{6}, C\sin \frac{5\pi }{6}) \\
\hspace{1.5in} = C(-\ttfrac{1}{2} +t, \ttfrac{\sqrt{3}}{2}),
\end{array} 
\end{displaymath}
which is a horizontal side of a regular hexagon centered at $(0,0)$ with side length $C$. Since $\gamma $ is an integral curve of $\frac{\partial }{\partial x}$ and hence is a geodesic on $(({\R }^2)^{\vvee}, \varepsilon )$, it takes time $C = \int^1_0 |\frac{\dee \gamma }{\dee t}|\, \dee t$ 
to go from $(C\cos \frac{4\pi}{6}, C\sin \frac{4\pi}{6} )$ to 
$(C \cos \frac{5\pi }{6}, C\sin \frac{5\pi }{6})$. Thus $\gamma $ runs off 
$({\R }^2)^{\vvee}$ in finite time. \hfill $\square $ \medskip

We remove the incompleteness of the geodesic vector field on $(K, {\gamma }_{|K})$ by imposing the 
following condition: when a geodesic starting at a point in 
$\mathrm{int} \, K$ meets $\partial K$ at a point on an open edge, it undergoes a reflection in that edge; otherwise 
it meets $K$ at a vertex, where it reverses its motion. Such motions in the stellated regular hexagon 
$K$, called \emph{billiard motions}, are defined for all time. 

\section{An affine model of $\boldsymbol{S^{\dagger}}$}

In this section we construct an affine model of the affine Riemann 
surface $S^{\dagger}$. In other words, we find a discrete subgroup 
$\mathfrak{G}$ of the $2$-dimensional Euclidean group $\mathrm{E}(2)$, which acts 
freely and properly on $\C \setminus {\mathbb{V}}^{+}$ and has 
the stellated regular hexagon $K$ less its vertices and center, as its fundamental domain. Here ${\mathbb{V}}^{+}$ is a discrete subset of $\C$ formed by translating the vertices of $K$ and its center. After constructing the identification space 
$(\C \setminus {\mathbb{V}}^{+})^{\sim }$, we show that its $\mathfrak{G}$ orbit space 
$(\C \setminus {\mathbb{V}}^{+})^{\sim }/\mathfrak{G}$ is holomorphically 
diffeomorphic to $S^{\dagger}$. \medskip 

First we specify the set ${\mathbb{V}}^{+}$. We label the edges of the regular hexagon $H$ of the stellated 
regular hexagon $K$ in order by $\{ 0, 1, \ldots , 5\} $. Reflecting the regular hexagon $H$ in its edge $k_0$ gives the regular hexagon $H_{k_0}$, which uniquely determines the reflected stellated regular hexagon $K_{k_0}$, see figure 4.

\par \noindent \hspace{.75in}\begin{tabular}{l}
\includegraphics[width=250pt]{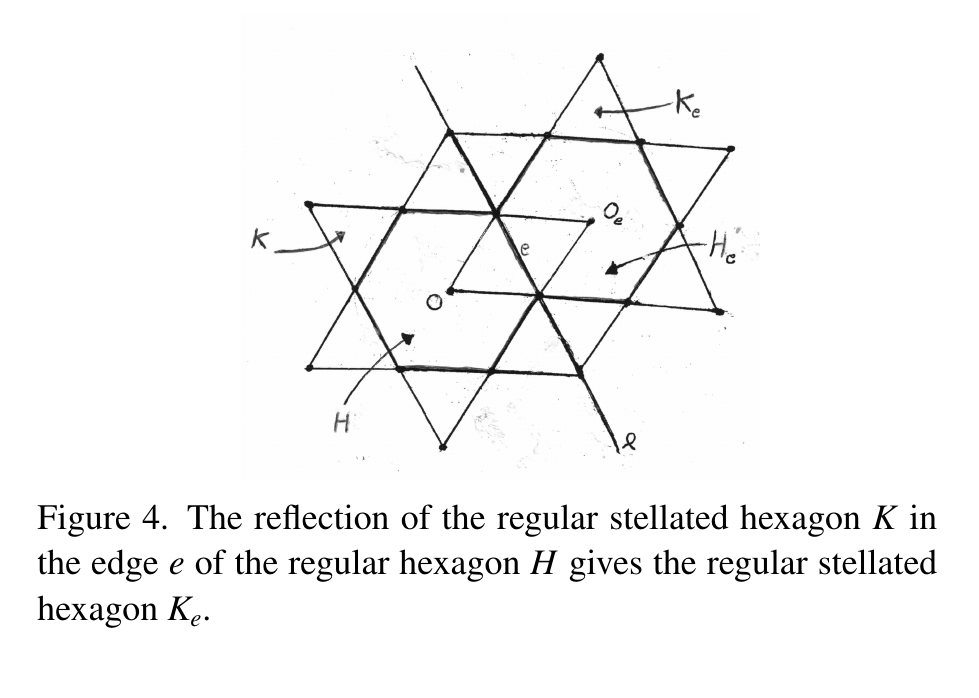}
\vspace{-.1in}
\end{tabular}  
\par \noindent Then we reflect in the edge $k_1$ of 
$K^{\ast}_{k_0}$ giving the stellated hexagon $K^{\ast }_{k_0k_1}$. After $n+1$ repetitions we get $K^{\ast }_{k_0k_1 \cdots k_n}$. Repeat this indefinitely and obtain the triangle tiling in figure 5. \medskip %

\vspace{-.1in}\noindent \hspace{.75in}\begin{tabular}{l}
\includegraphics[width=250pt]{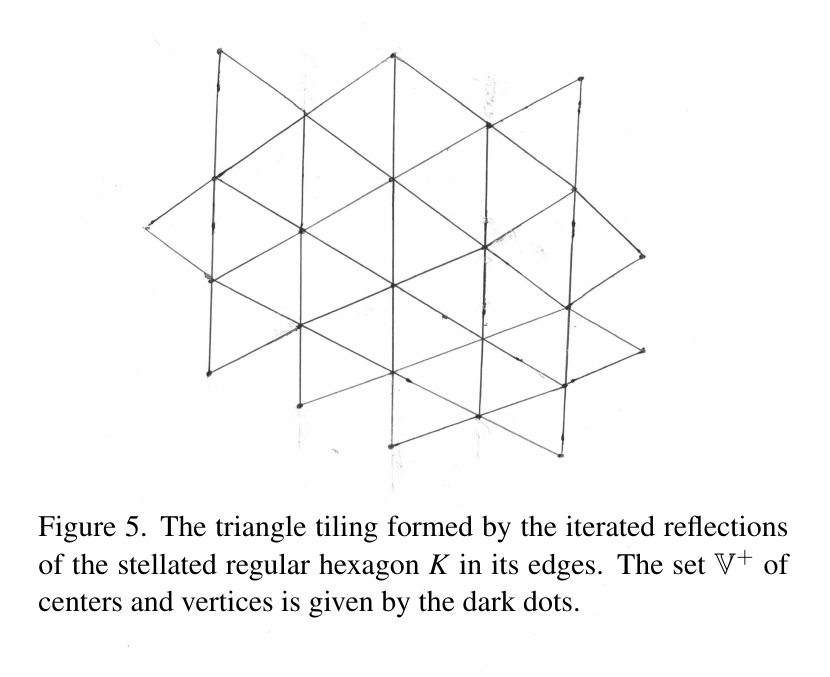}
\end{tabular} 

Consider the translations 
\begin{equation}
{\tau }_k: \C \rightarrow \C : z \mapsto z + u_k,
\label{eq-5onea}
\end{equation}%
where $u_k= \sqrt{3}CR^k$ for $k= 0,1, \ldots , 5$. \medskip 

The set of centers of the iteratively reflected stellated regular hexagons is 
\begin{displaymath}
\{ ({\tau }^{{\ell }_0}_0 \comp \cdots \comp {\tau }^{{\ell }_5}_5)(0) \in \C \setrule \, ({\ell }_0, {\ell }_1 , \ldots , {\ell }_5) \in {\Z }^6 \}. 
\end{displaymath}
The set of vertices of $K$ is 
\begin{displaymath}
V = \{ V_{2\ell}= CR^{2\ell}, \, \,  V_{2\ell -1} = \ttfrac{\sqrt{3}}{2}CR^{2\ell -1}, \, \, \mbox{for $\ell =0,1, \ldots , 5$} \} . 
\end{displaymath}
The set of centers and vertices the hexagonal tiling of $\C$ is 
\begin{displaymath}
{\mathbb{V}}^{+} = \{ {\tau }^{{\ell }_0}_0 \comp \cdots \comp {\tau }^{{\ell }_5}_5(V^{+}) 
\in \C \setrule \, V^{+} = V \cup \{ 0 \} \, \& \, ({\ell }_0, \ldots , {\ell }_5) \in {\Z }^6 \}   ,
\end{displaymath}
see figure 5. \medskip

Let $\mathcal{T}$ be the abelian subgroup of $\mathrm{E}(2)$ generated by 
the translations ${\tau }_k$, $k=0, \ldots 5$, see (\ref{eq-5onea}). $\mathcal{T}$ is 
isomorphic to the additive group $\{ \sum^5_{j=0}{\ell }_ju_j \in \C 
\setrule \, ({\ell }_1, $ $\ldots , {\ell }_6) $ $\in {\Z }^6 \}$. By definition 
$K^{\ast } \setminus O$ is a fundamental domain for the action of $\mathcal{T}$ 
on $\C \setminus {\mathbb{V}}^{+}$. Let $\mathfrak{G} = \widetilde{G} \ltimes 
\mathcal{T} \subseteq \widetilde{G} \times \mathcal{T}$. An element of 
$\mathfrak{G}$ is the affine mapping 
\begin{displaymath}
(R^j, u_k) : \C \rightarrow \C: z \mapsto R^jz + u_k, \quad 
\mbox{for $j,k \in \{0, \ldots , 5\}$.} 
\end{displaymath}
Multiplication $\cdot $ in $\mathfrak{G}$ is given by 
\begin{align*}
(R^h, u_{\ell })\cdot (R^j, u_k) & = (R^{(h+j)\! \bmod 6}, R^hu_k + u_{\ell }) = 
(R^{(h+j)\! \bmod 6}, u_{(k + \ell ) \bmod 6}+u_k) . 
\end{align*}
The group $\mathfrak{G}$ acts on $\C $ just as $\mathrm{E}(2)$ does, namely, by 
affine orthogonal mappings. Denote this action by 
\begin{equation}
\psi : \mathfrak{G} \times \C \rightarrow \C : \big( (g,\tau ), z \big) 
\mapsto \tau \big( g(z) \big) . 
\label{eq-5two}
\end{equation}

\noindent \textbf{Lemma 5.1} The set ${\mathbb{V}}^{+}$ is invariant under the 
$\mathfrak{G}$ action $\psi $ (\ref{eq-5two}). \medskip 

\noindent \textbf{Proof.} Let $v \in {\mathbb{V}}^{+}$. For some 
$(k_0, \ldots , k_n) \in {\{ 0, \ldots, 5 \} }^{n+1}$ and some $w \in V^{+}$ 
\begin{displaymath}
v ={\tau }_{k_n} \comp \cdots \comp {\tau }_{k_0}(w) = 
{\psi }_{(e,u)}(w) , 
\end{displaymath}
where $u = \sum^5_{j=0} {\ell }_j u_j$ with $({\ell }_0, \ldots , {\ell }_5) \in {\Z }^6$. 
Each ${\ell }_j$ is uniquely determined by $(k_0, \ldots , k_n) \in 
{\{ 0, \ldots , 5 \}}^{n+1}$.  For $(R^j, u') \in \mathfrak{G}$ with 
$j =0, \ldots , 5$ and $u' \in \mathcal{T}$ 
\begin{align*}
{\psi }_{(R^j,u')}(v) & = {\psi }_{(R^j,u')}\comp {\psi }_{(e, u)}(w) = 
{\psi }_{(R^j,u') \cdot (e, u)}(w)  \\
& = {\psi}_{(R^j, R^j u +u')}(w) = {\psi }_{(e, R^ju+u')\cdot (R^j,0)}(w) \\
& = {\psi }_{(e, R^ju+u')}\big( {\psi }_{(R^j,0)}(w) \big) = 
{\psi }_{(e, R^ju +u')}(w'),
\end{align*}
where $w' = {\psi }_{(R^j,0)}(w) \in V^{+}$. Since $(e,u') \in \mathcal{T}$ and 
\begin{displaymath}
R^ju = \sum^5_{k=0} {\ell }_k R^ju_k = 
\sum^5_{k=0} {\ell}_ku_{(k+j) \bmod 6}, 
\end{displaymath}
it follows that $(e, R^ju +u') \in \mathcal{T}$. Hence 
${\psi }_{(R^j ,u)}(v) \in {\mathbb{V}}^{+}$, as desired. \hfill $\square $ \medskip  

\noindent \textbf{Lemma 5.2} The action of $\mathfrak{G}$ on $\C \setminus {\mathbb{V}}^{+}$ is free. \medskip 

\noindent \textbf{Proof.} Suppose that for $v \in \C \setminus {\mathbb{V}}^{+}$ 
there is $(R^j,u) \in \mathfrak{G}$ such that $v = {\psi }_{(R^j,u)}(v)$. Then 
$v$ lies in the stellated hexagon $K^{\ast }_{k_0 \cdots k_n}$ given 
by reflecting $K^{\ast }$ in the successive edges $k_0$, $k_1$, \ldots , 
$k_{n}$. So for some $v' \in K^{\ast }$ 
\begin{displaymath}
v = ({\tau }_{k_n} \comp \cdots \comp {\tau }_{k_0})(v') = {\psi }_{(e, u')}(v'),  
\end{displaymath}
where $u' =\sum^5_{j=0}{\ell }_j u_j$ for some $({\ell }_0, \ldots , {\ell }_5) \in 
{\Z }^6$, which is uniquely determined by $(k_0, \ldots , k_n) \in 
{\{ 0, \ldots , 5 \}}^{n+1}$. Thus 
\begin{displaymath}
{\psi }_{(e,u')}(v') = {\psi}_{(R^j,u) \cdot (e, u')}(v') = {\psi }_{(R^ju, R^ju'+u)}(v'), 
\end{displaymath}
which implies $R^j = e$, that is, $j =0 \bmod 6$ and $u' = R^j u' + u = u'+u$, 
that is, $u =0$. Hence $(R^j,u) = (e,0)$, which is the identity element of 
$\mathfrak{G}$. \hfill $\square $ \medskip 

\noindent \textbf{Lemma 5.3} The action of $\mathcal{T}$ (and hence 
of $\mathfrak{G}$) on 
\begin{displaymath}
\C \setminus 
{\mathbb{V}}^{+} = K^{\ast }\bigcup \big( \cup_{n \ge 0}\cup_{0 \le j \le n } 
\cup_{0 \le k_j \le 5} K^{\ast }_{k_0 \cdots k_n} \big)
\end{displaymath} 
is transitive. \medskip 

\noindent \textbf{Proof.} Let $K^{\ast }_{k_0 \cdots k_n}$ and 
$K^{\ast }_{k'_0 \cdots k'_{n'}}$ lie in $\C \setminus {\mathbb{V}}^{+}$. 
Since $K^{\ast }_{k_0 \cdots k_n} = ({\tau }_{k_n} \comp \cdots 
\comp {\tau }_{k_0})(K^{\ast })$ and $K^{\ast }_{k'_0 \cdots k'_{n'}} = 
({\tau }_{k'_{n'}} \comp \cdots \comp {\tau }_{k'_0})(K^{\ast })$, it follows that 
\begin{displaymath}
K^{\ast }_{k'_0 \cdots k'_{n'} } = ({\tau }_{k'_{n'}} \comp \cdots \comp {\tau }_{k'_0}) 
\comp ( {\tau }_{k_n} \comp \cdots \comp {\tau }_{k_0})^{-1}(K^{\ast }) . 
\end{displaymath}
Thus the action of $\mathcal{T}$ on $\C \setminus {\mathbb{V}}^{+}$ is 
transitive. \hfill $\square $ \medskip 

The action of $\mathfrak{G}$ on $\C \setminus {\mathbb{V}}^{+}$ 
is proper because $\mathfrak{G}$ is a discrete subgroup of $\mathrm{E}(2)$, 
having no accumulation points. \medskip 

Let $E_{k_0 \cdots k_n} = ({\tau }_{k_n} \comp \cdots \comp {\tau }_{k_0})(E) 
\in K^{\ast }_{k_0 \cdots k_n}$. Then $E_{k_0 \cdots k_n}$ 
is an \emph{open edge} in $\C \setminus {\mathbb{V}}^{+}$ of the stellated hexagon $K^{\ast }_{k_0 \cdots k_n}$. 
So $\mathfrak{E} = \{ E_{k_0 \cdots k_n} \setrule \, 
n \ge 0, \, \, 0 \le j \le n \, \, \& \, \, 0 \le k_j \le 5 \} $ be the set of open edges of 
$\C \setminus {\mathbb{V}}^{+}$. By definition $\mathfrak{E}$ is invariant under translations in $\mathcal{T}$. \medskip 

\noindent \textbf{Lemma 5.4} The $\mathfrak{G}$ action $\psi $ (\ref{eq-5two}) leaves the subset $\mathfrak{E}$ of $\C \setminus {\mathbb{V}}^{+}$ invariant. \medskip 

\noindent \textbf{Proof.} Let $F \in \mathfrak{E}$. For some 
$(k_0, \ldots , k_n) \in {\{ 0, \ldots , 5 \}}^{n+1}$ and some open edge 
$E$ of $K^{\ast}$ 
\begin{displaymath}
F =({\tau }_{k_n} \comp \cdots \comp {\tau }_{k_0}(E) = 
{\psi }_{(e,u')}(E) , 
\end{displaymath}
where $u' = \sum^5_{j=0} {\ell}_j u_j$ with $({\ell}_0, \ldots , {\ell }_5) \in {\Z }^6$. 
Each ${\ell }_j$ is uniquely determined by $(k_0, \ldots , k_n) \in 
{\{ 0, \ldots , 5 \} }^{n+1}$.  For 
$(R^j, u) \in \mathfrak{G}$ with $j =0,1, \ldots , 5$ and $u \in \mathcal{T}$ 
\begin{align*}
{\psi }_{(R^j,u)}(F) & = {\psi }_{(R^j,u)}\comp {\psi }_{(e, u')}(E) = 
{\psi }_{(R^j,u) \cdot (e, u')}(E)  \\
& = {\psi}_{(R^j, R^j u' +u)}(E) = {\psi }_{(e, R^ju'+u)\cdot (R^j,0)}(E) \\
& = {\psi }_{(e, R^ju'+u)}\big( {\psi }_{(R^j,0)}(E) \big) = 
{\psi }_{(e, R^ju' +u)}(E'),
\end{align*}
where $E' = {\psi }_{(R^j,0)}(E) = R^jE$ is an open edge of $K^{\ast }$. 
Since $(e,u) \in \mathfrak{G}$ and 
\begin{displaymath}
R^ju' = \sum^5_{i=0} {\ell}_i R^ju_i = 
\sum^5_{i=0} {\ell}_i u_{(j+i) \bmod 6}, 
\end{displaymath}
it follows that $(e, R^ju' +u) \in \mathcal{T}$. Hence 
${\psi }_{(e, R^ju'+u)}(E') \in \mathfrak{E}$. So ${\psi }_{(R^j,u)}(F) \in 
\mathfrak{E}$, as desired. \hfill $\square $ \medskip  

\noindent We say that two open edges in $\mathfrak{E}$ are \emph{equivalent} if they are parallel in $\C$. The $\mathfrak{G}$ action $\psi $ (\ref{eq-5two}) on 
$\mathfrak{E}$ preserves the relation of equivalence of edges. \medskip 

\noindent \textbf{Proposition 5.5} The identification space 
$(\C \times {\mathbb{V}}^{+})^{\sim }$ formed by identifying equivalent open 
edges in $\mathfrak{E}$ is equal to the identification space 
$(K^{\ast } \setminus O)^{\sim}$ formed by identifying equivalent open edges 
of the stellated regular hexagon $K^{\ast }$ less its vertices. \medskip 

\noindent \textbf{Proof.} This follows from the observation that 
every equivalence class of open edges in $\mathfrak{E}$ contains a 
unique equivalence class of open edges of $K^{\ast }$. \hfill $\square $ \medskip 

\noindent \textbf{Corollary 5.5A} The $\mathfrak{G}$ orbit space 
$(\C \times {\mathbb{V}}^{+})^{\sim }/\mathfrak{G}$ is equal to 
the $\widetilde{G}$ orbit space $(K^{\ast } \setminus O)^{\sim }/\widetilde{G}$, which is a $1$-dimensional complex manifold. \medskip 

\noindent \textbf{Proof.} Since the action of $\mathfrak{G}$ on 
$\C \times {\mathbb{V}}^{+}$ is proper and free, it induces a free and 
proper action on the identification space $(\C \times {\mathbb{V}}^{+})^{\sim }$. 
Hence the orbit space $(\C \times {\mathbb{V}}^{+})^{\sim }/\mathfrak{G}$ is 
a complex manifold, which is equal to the orbit space 
$(K^{\ast } \setminus O)^{\sim }/\widetilde{G}$. \hfill $\square $ \medskip 

The orbit space $(\C \times {\mathbb{V}}^{+})^{\sim }/\mathfrak{G} = 
{\widetilde{S}}^{\dagger}$ is 
an affine model of the affine Riemann surface $S^{\dagger}$, since 
$(K^{\ast } \setminus O)^{\sim }/\widetilde{G}$ is holomorphically diffeomorphic 
to $S^{\dagger}$. \medskip 
 
\noindent \textbf{Theorem 5.6} The image of a $\widetilde{\mathcal{G}}$ invariant geodesic on $(S^{\dagger}, \Gamma )$ under the  
developing map $\delta $ (\ref{eq-sixstar}) is a $\widetilde{G}$ invariant 
billiard motion in $K^{\ast }$. \medskip 

\noindent \textbf{Proof.} Because ${\mathcal{R}}^j$ and 
$R^j$ are isometries of $(S^{\dagger}, \Gamma )$ and 
$(K^{\ast }, {\gamma }_{| K^{\ast }})$, respectively, it follows that the surjective map 
$\delta : (S^{\dagger}, \Gamma ) 
\rightarrow (K^{\ast }, {\gamma }_{|K^{\ast }} )$ is an isometry, because 
$\delta $ pulls back the metric ${\gamma }_{|K^{\ast }}$ to the metric 
$\Gamma $.  Hence $\delta $ is a local developing map. Using the local inverse of 
$\delta $, it follows that a billiard motion in $\mathrm{int}(K^{\ast } \setminus  \mathrm{O}  )$ is mapped onto a geodesic in 
$(S^{\dagger}, \Gamma )$, which is possibly broken at the 
points $({\xi }_i, {\eta }_i) = {\delta }^{-1}_{K^{\ast}}(p_i)$. Here $p_i \in \partial K^{\ast}$ are the points where the billiard motion undergoes a reflection or a reversal. But the geodesic on $S^{\dagger}$ is smooth at $({\xi }_i, {\eta }_i)$ since the geodesic vector field is holomorphic on $S^{\dagger}$. Thus the image of a geodesic under the developing map $\delta $ is a billiard motion. If the geodesic is 
$\widetilde{\mathcal{G}}$ invariant, then the image billiard motion is $\widetilde{G}$ invariant, because the developing map intertwines the $\widetilde{\mathcal{G}}$ action on $S^{\dagger}$ with the $\widetilde{G}$ action on $K^{\ast }$. \hfill 
$\square $ \medskip 

\par\noindent \hspace{.25in}\begin{tabular}{l}
\vspace{-.2in}\\
\includegraphics[width=300pt]{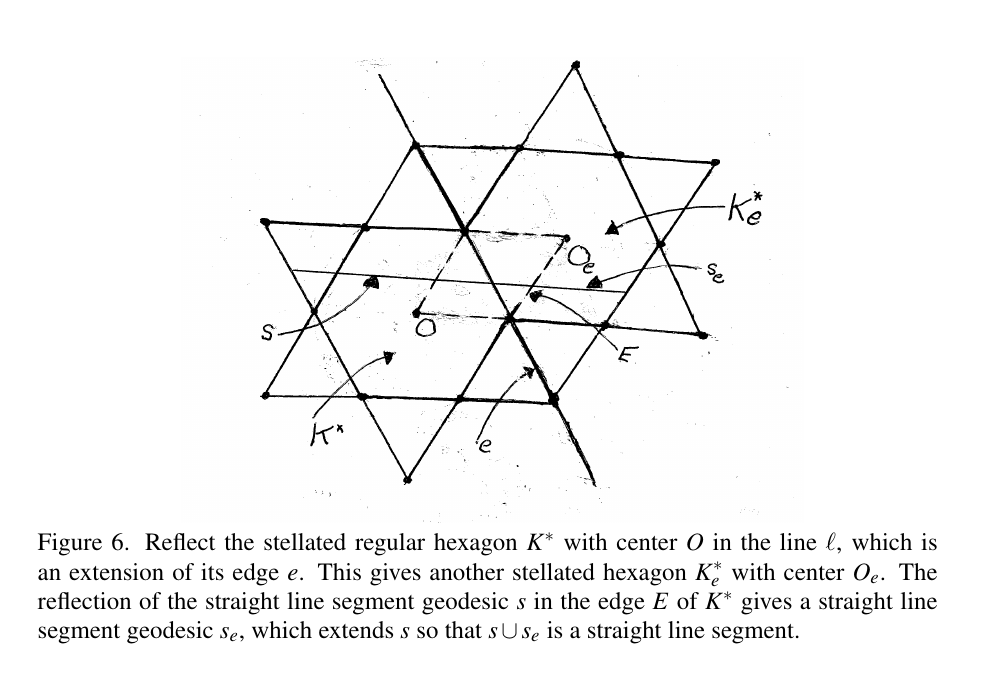}
\vspace{-.1in}
\end{tabular}

\noindent \textbf{Theorem 5.7} Under the restriction of the mapping 
\begin{equation}
\nu = \sigma \comp \rho : \C \setminus {\mathbb{V}}^{+}   \rightarrow  (\C \setminus {\mathbb{V}}^{+} )^{\sim }/\mathfrak{G} = {\widetilde{S}}^{\dagger}
\label{eq-5three}
\end{equation} 
to $K^{\ast }\setminus  \mathrm{O}  $ the image of a $\widetilde{G}$ invariant 
billiard motion ${\lambda }_z$ is a smooth geodesic ${\widehat{\lambda }}_{\nu (z)}$ on $({\widetilde{S}}^{\dagger} , \widehat{\gamma } ) $. Here 
$\widehat{\gamma }$ is the metric on ${\widetilde{S}}^{\dagger}$ whose pull back under the mapping $\nu $ (\ref{eq-5three}) is the metric ${\gamma }_{|K^{\ast }}$. \medskip 

\noindent \textbf{Proof.} Since the Riemannian metric 
$\gamma $ on $\C $ is invariant under the group of Euclidean motions, the Riemannian metric ${\gamma }|_{K^{\ast } \setminus  \mathrm{O} }$ on $K^{\ast } \setminus \mathrm{O}$ is $\widetilde{G}$-invariant. Hence 
${\gamma }_{|K^{\ast } \setminus  \mathrm{O}  }$ is invariant 
under the reflection $S_m$ for $m \in \{ 0,1, \ldots , 5 \} $. So 
${\gamma }|_{K^{\ast } \setminus  \mathrm{O}  }$ piece together to give a Riemannian metric ${\gamma }^{\sim }$ on the identification space 
$(K^{\ast } \setminus  \mathrm{O}  )^{\sim}$. In other words, the pull back of 
${\gamma }^{\sim }$ under the map ${\rho }|_{K^{\ast } \setminus \mathrm{O}  }: 
K^{\ast } \setminus  \mathrm{O}  \rightarrow (K^{\ast } \setminus 
 \mathrm{O} )^{\sim }$, which identifies equivalent edges of $K^{\ast }$, is the metric 
${\gamma }|_{K^{\ast } \setminus  \mathrm{O}  }$. Since ${\rho }|_{K^{\ast } \setminus 
\mathrm{O}  }$ intertwines the $\widetilde{G}$-action on 
$K^{\ast } \setminus  \mathrm{O} $ with the $\widetilde{G}$-action on 
$(K^{\ast } \setminus \mathrm{O} )^{\sim }$, the metric ${\gamma }^{\sim }$ is 
$\widetilde{G}$-invariant. It is flat because the metric $\gamma $ is flat. So 
${\gamma }^{\sim }$ induces a flat Riemannian metric $\widehat{\gamma }$ 
on the orbit space $(K^{\ast } \setminus  \mathrm{O}  )^{\sim } /\widetilde{G} = 
{\widetilde{S}}^{\dagger}$. Since the billiard motion ${\lambda}_z$ is 
a $\widetilde{G}$-invariant broken geodesic on $(K^{\ast } \setminus  \mathrm{O} , 
{\gamma }_{K^{\ast } \setminus  \mathrm{O}  })$, it gives rise to a \emph{continuous} broken geodesic ${\lambda }^{\sim }_{\rho (z)}$ on $((K^{\ast } \setminus 
\mathrm{O} )^{\sim }, {\gamma }^{\sim })$, which is $\widetilde{G}$-invariant. Thus 
${\widehat{\lambda }}_{\nu (z) } = \nu ({\lambda}_z)$ is a piecewise smooth geodesic on the smooth $\widetilde{G}$-orbit space 
$((K^{\ast } \setminus  \mathrm{O} )^{\sim }/\widetilde{G} 
= {\widetilde{S}}^{\dagger}, \widehat{\gamma })$. We need only show that ${\widehat{\lambda }}_{\nu (z)}$ is smooth. To see this we argue as follows. Let $s \subseteq K^{\ast }$ be a closed segment of a billiard motion ${\gamma }_z$, that does not meet a vertex of 
$K$. Then we may assume that $s$ is a horizontal straight line motion in $K$. Suppose that $E_{k_0}$ is the edge of $K^{\ast }$, perpendicular to the direction $u_{k_0}$, which is first met by $s$ and let 
$P_{k_0}$ be the meeting point. Let $S_{k_0}$ be the reflection in $E_{k_0}$. 
The continuation of the motion $s$ at $P_{k_0}$ is the horizontal line $RS_{k_0}(s)$ in $K^{\ast }_{k_0}$. Recall that $K^{\ast }_{k_0}$ is the translation 
of $K^{\ast }$ by ${\tau }_{k_0}$. Using a suitable sequence of reflections in the edges of a suitable $K^{\ast }_{k_0 \cdots k_n}$ each followed by a rotation $R$ and then a translation in $\mathcal{T}$ corresponding to their origins, we extend $s$ to a smooth straight line $\lambda $ in $\C \setminus {\mathbb{V}}^{+}$, see figure 6. The line 
$\lambda $ is a geodesic in $(\C \setminus {\mathbb{V}}^{+}, \gamma |_{\C \setminus {\mathbb{V}}^{+}})$, which in $K^{\ast }$ has image 
${\widehat{\lambda }}_{\nu (z)}$ under the $\mathfrak{G}$-orbit map $\nu $ (\ref{eq-5three}) that is a smooth geodesic on $({\widetilde{S}}^{\dagger}, \widehat{\gamma } ) $. The geodesic $\nu (\lambda )$ starts at $\nu (z)$. Thus the smooth geodesic $\nu (\lambda )$ and the geodesic ${\widehat{\lambda}}_{\nu (z)}$ are equal. In other words, ${\widehat{\lambda}}_{\nu (z)}$ is a smooth geodesic. \hfill $\square $ \medskip 

The affine orbit space ${\widetilde{S}}^{\dagger} = 
(\C \setminus {\mathbb{V}}^{+})^{\sim} /\mathfrak{G}$ with flat Riemannian metric 
$\widehat{\gamma }$ is the \emph{affine} analogue of the Poincar\'{e} model of 
the Riemann surface $S^{\dagger}$ as an orbit space of a discrete subgroup of 
$\mathrm{PGl}(2, \C)$ acting on the unit disk in $\C $ with the Poincar\'{e} metric, 
see \cite{weyl}.

\end{document}